\documentclass[12pt]{article}
\usepackage[final]{graphics}
\usepackage{amsthm,amsfonts,epsfig}
 \usepackage{latexsym, amssymb}

\def\ds{\displaystyle}
\def\Lra{\Leftrightarrow}
\def\ov{\over}\def\pl{\partial}
\def\smoc{{\cal C}^\infty_c}
\def\td{\tilde}\def\la{\langle}\def\ra{\rangle}
\def\ep{\epsilon}
\def\Gam{\Gamma}\def\gam{\gamma}
\def\shp{\sharp}
\def\msk{\medskip}
\newtheorem{theorem}{Theorem}[section]

\newtheorem{lemma}{Lemma}[section]
\newtheorem{cor}{Corollary}[section]

\title{Hypergeometric functions and the Tricomi operator} 
\author{J.~Barros-Neto\footnote{Partially supported by NSF, Grant \# INT 0124940}\\
Rutgers University, Hill Center\\
110 Frelinghuysen Rd, Piscataway, NJ 08854-8019\\
e-mail: jbn@math.rutgers.edu
\and
Fernando Cardoso\footnote{Partially supported by CNPq (Brazil)}\\
Departmento de Matem\'atica, Universidade Federal de Pernambuco\\
50540-740 Recife, Pe, Brazil\\
e-mail: fernando@dmat.ufpe.br} 
\date{} 
\begin{document}  
\maketitle
\begin{abstract} In this paper we show how certain hypergeometric functions play an important role in finding
fundamental solutions for a generalized Tricomi operator.
\end{abstract}

\section{Introduction}\setcounter{equation}{0}\label{int}
In this article we consider the operator 
\begin{equation}\label{gn1}
{\cal T}=y\Delta_x +{\pl^2\over\pl y^2},
\end{equation}
in $\mathbb R^{n+1},$ where $\Delta_x=\sum_{j=1}^n\ds{\pl^2\ov\pl x_j^2}, $ $n\geq1.$ This is a natural
generalization of the classical Tricomi operator in $\mathbb R^2$ already considered by us in the article
\cite{bc} where it was called {\em generalized Tricomi operator.} 

In that article we obtained, by the method of partial Fourier transformation, explict expressions for fundamental
solutions to $\cal T,$ relative to points on the hyperplane $y=0.$ That lead us to calculate inverse Fourier
transforms of Bessel functions which, in turn, revealed the importance of certain hypergeometric functions
(depending on the ``space dimension'' $n$) that are intimately related to the operator $\cal T.$

In the present article we look for fundamental solutions of $\cal T$ relative to an arbitrary point $(x_0,y_0),$  
located in the hyperbolic region ($y<0$) of the operator, and which are supported by the ``forward''   
characteristic conoid of $\cal T$ with vertex at $(x_0,y_0).$ We follow the method of S.~Delache and J.~Leray   
in \cite{dl} where they introduced {\em hypergeometric distributions}, a notion also considered
by I.~M.~Gelfand and G.~E.~Shilov in \cite{gs}.

The plan of this article is the following.  In Section \ref{gn} we deal with preliminary material that is needed
throughout the paper. Hypergeometric distributions are introduced in Section \ref{hd} where we obtain the basic 
formula (\ref{hd14}) which is used in Sections \ref{cl}, \ref{ev}, and \ref{od} to obtain fundamental solutions
respectively in the cases $n=1$ (the classical Tricomi operator), $n$ even, and $n$ odd.  The case $n$ odd 
$\geq 3$ differs from the other two cases by the fact that the fundamental solution is then a sum of two terms one
supported by the ``forward'' conoid (as in the cases $n=1$ and $n$ even) and another supported by the boundary
of the conoid. In Section \ref{cl} we also show how to derive from the method used in this paper the results obtained
previously by Barros-Neto and Gelfand in \cite{bg1}. Finally, in the Appendix we prove or indicate the proof of
results mentioned and utilized in Section \ref{cl}. 

\section{Preliminaries}\setcounter{equation}{0}\label{gn}
Let $\cal T$ be the operator given by (\ref{gn1}) at the begining of Section \ref{int}. Defining the modified
gradient 
\begin{equation}\label{gn1a}
\dot\nabla u=\la yu_{x_1},\ldots,yu_{x_n},u_y\ra,
\end{equation}
one verifies that ${\cal T}u=\mbox{div}(\dot\nabla u)$ and that
\begin{equation}\label{gn1b}
\int\!\!\int_D(u{\cal T}v-v{\cal T}u)\,dV=\int_{\pl D}(u\dot\nabla v-v\dot\nabla u)\cdot\vec n\,dS
\end{equation}
for all smooth $u$ and $v$ on the closure of an open bounded domain $D$ with smooth boundary      
$\pl D.$

If $y<0$ and we set $t=2(-y)^{3/2}/3>0,$ then the change of variables
\begin{equation}\label{gn2}
x=x,\ t=2(-y)^{3/2}/3\quad \Lra\quad x=x,\ 
 y=-({2/3})^{-2/3}t^{2/3},
\end{equation}
whose Jacobian is
\begin{equation}\label{gn3}
{\pl(x,y)\ov\pl(x,t)}=-({2/3})^{1/3}t^{-1/3},
\end{equation}
transforms ${\cal T}$ into the operator 
\begin{equation}\label{gn4}
2({3t\ov2})^{2/3}\,{\cal T}_h,
\end{equation}
with
\begin{equation}\label{gn5}
{\cal T}_h={1\ov2}({\pl^2\ov\pl t^2}-\Delta_x) +{1\ov 6t}{\pl\ov\pl t}.
\end{equation}
We call ${\cal T}_h$ the {\em reduced hyperbolic} Tricomi operator. Its formal adjoint is
\begin{equation}\label{gn6}
{\cal T}^*_h={1\ov2}({\pl^2\ov\pl t^2}-\Delta_x) -{1\ov 6t}{\pl\ov\pl t} +{1\ov6t^2}.
\end{equation}

It is a matter of verification that
$$
{\cal T}^*_h(t^{1/3}u)=t^{1/3}{\cal T}_h(u).
$$
Thus, if $u$ is a solution of ${\cal T}_h(u)=0,$ then $v=t^{1/3}u$ is a solution of ${\cal T}^*_hv=0,$ and
conversely. Moreover, suppose that $E(x,t;0,t_0),$ with $t_0\neq0,$ is a fundamental solution of ${\cal T}_h$
relative to the  point $(0,t_0),$ that is,
$$
{\cal T}_hE=\delta(x,t-t_0),
$$
then  $(t/t_0)^{1/3}E$ is a fundamental solution of  ${\cal T}^*_h$ relative to the same point.

We now recall the definition of the distribution $\chi_q(s)$ (see \cite{dl, gs,sch}).  Let $q\in\mathbb C$ be
such that $\mbox{Re}\, q>-1.$  The locally integrable function  
\begin{equation}\label{gn7}
\chi_q(s)={s^q\ov\Gamma(q+1)}\ \mbox{if}\ s>0,\quad \chi_q(s)=0\ \mbox{if}\  s\leq 0
\end{equation}
defines a distribution in $\mathbb R$ that depends analytically on $q$ and that extends by analytic continuation
to an entire function of $q.$ We have
$$
\chi_q(s)={d\ov ds}\chi_{q+1}(s).
$$
Moreover, $\chi_q(s)$ is positive homogeneous of degree $q$ and Euler's formula holds
\begin{equation}\label{gn8}
s\chi_{q-1}(s)=q\chi_q(s).
\end{equation}
We also have that $\chi_q(s)=\delta^{(-q-1)}(s)$ if $q$ is an integer $<0$ (see \cite{gs}). 

Consider the function 
\begin{equation}\label{gn9}
k(x,t-t_0) = \left\{\begin{array}{ll}
         (t-t_0)^2-|x|^2 & \mbox{if $t-t_0>|x|$}\\
          \\
         0 & \mbox{if $t-t_0\leq |x|$},
\end{array} \right.
\end{equation}
defined in the whole of $\mathbb R^{n+1}.$ Since $k(x,t-t_0)$ is positive in the semi-cone
$C=\{(x,t)\in\mathbb R^{n+1}:t-t_0>|x|$\}, and identically zero outside of $C,$ it follows that
$\chi_q(k(x,t-t_0))$ (which, for simplicity and when no confusion is possible, we denote by $\chi_q(k(\cdot))$) 
is a distribution in $\mathbb R^{n+1}$ which is an entire analytic function of $q\in\mathbb C.$ In particular, if    
$q$ is an integer $<0,$ then
$$
\chi_q(k(\cdot))=\delta^{(-q-1)}((k(\cdot)))
$$
is a distribution concentrated on the boundary of $C$ (see \cite{gs}).

\section{Hypergeometric distributions}\setcounter{equation}{0}\label{hd}

Our aim is to find fundamental solutions for the Tricomi operator ${\cal T}$ relative to an arbitrary point 
$(0,b), b<0,$ in $\mathbb R^{n+1},$ that is, a distribution $E(x,y;0,b)$ defined in $\mathbb R^{n+1}$ so that 
${\cal T}E=\delta(x,y-b).$ In guise of motivation, suppose that $E(x,y;0,b)$ is locally integrable function. Then in
view of formulas (\ref{gn2}), (\ref{gn3}), (\ref{gn4}),  and (\ref{gn5}), we may write
\begin{eqnarray}\label{hd1}
\phi(0,b)&=&\la{\cal T}E,\phi\ra=\int_{{\mathbb R}^{n+1}}E(x,y;0,b){\cal T}\phi\,dx\,dy\\
&=& 2({3\ov2})^{1/3}\int_{{\mathbb R}^{n+1}}t^{1/3}E^\sharp(x,t;0,t_0){\cal T}_h
          \psi(x,t)\,dx\,dt,\nonumber
\end{eqnarray}
where in the last formula $E^\sharp(x,t;0,t_0)$ and $\psi(x,t)$ denote, respectively, $E(x,y;0,b)$ and $\phi$ in
the  variables $x$ and $t,$ and we have set $t_0=2(-b)^{3/2}/3.$ Thus, our problem reduces to finding
fundamental  solutions relative to $(0,t_0)$ for the adjoint operator ${\cal T}^*_h,$ which according to our remark
in Section 1  is equivalent to finding fundamental solutions for ${\cal T}_h$ relative to the same point. 

The operator ${\cal T}_h$ belongs to a class of operators, called Euler--Poisson--Darboux operators,
studied by Delache and Leray in \cite{dl}, where they obtained explicit formulas for fundamental
solutions of those operators. For sake of completeness, we outline Delache and Leray's method in \cite{dl} relative
to the reduced hyperbolic Tricomi operator ${\cal T}_h,$ or more generally, the operator
\begin{equation}\label{hd2}
{\cal P}_\alpha={1\ov2}({\pl^2\ov\pl t^2}-\Delta_x) +{\alpha\ov t}{\pl\ov\pl t},
\end{equation}
where $\alpha\in\mathbb C.$ Note that ${\cal P}_\alpha$ remains invariant under the action of the group that
leaves $t$ unchanged and transforms $(x_1,\ldots,x_n)$ by translations. Since ${\cal P}_\alpha$ is homogeneous
of degree $-2$ and $\delta(x,t-t_0)$ is homogeneous of degree $-(n+1),$ a fundamental solution $E_\alpha$ of
${\cal P}_\alpha$ should be homogeneous  of degree $1-n.$  Monomials of the type
$$
t_0^{\alpha-j}t^{-\alpha-j}\chi_{j+1/2-n/2}(k(x,t-t_0)),
$$
have the desired homogeneity degree. On the other hand, if $\Box_{(x,t)}=({\pl^2/\pl t^2}-\Delta_x)$ denotes the
wave operator in $\mathbb R^{n+1},$ it is shown in \cite{dl} that 
\begin{equation}\label{hd3}
\Box_{(x,t)}({1\ov2}\pi^{1/2-n/2}\chi_{1/2-n/2}(k(x,t))=\delta(x,t),
\end{equation}
in other words, the distribution $\pi^{1/2-n/2}\chi_{1/2-n/2}(k(x,t))/2$ is a fundamental solution of the wave
operator.

As a consequence, it is natural to look for a fundamental solution to ${\cal P}_\alpha$ as a formal series
\begin{equation}\label{hd4}
E_\alpha(x,t;0,t_0)=\pi^{1/2-n/2}({t_0\ov t})^\alpha
   \sum_{j=0}^\infty c_j(t_0t)^{-j}\chi_{j+1/2-n/2}(k(\cdot)),
\end{equation}
with a suitable choice of the coefficients $c_j.$ By applying ${\cal P}_\alpha$ to both sides of (\ref{hd4}) one 
obtains, after routine calculations where the two identities
$$
\Box_{(x,t)}\chi_{j+1/2-n/2}(k(\cdot))=4j\chi_{j-1/2-n/2}(k(\cdot))
$$
and
$$
\ds{\pl\ov\pl t}[\chi_{j+1/2-n/2}(k(\cdot))]=2(t-t_0)\chi_{j-1/2-n/2}(k(\cdot))
$$
are used, the following result:
\begin{equation}\label{hd5}
{\cal P}_\alpha E_\alpha = c_0\delta(x,t-t_0) +
\end{equation}
$$
+\pi^{1/2-n/2}\sum_{j=1}^\infty\{{1\ov2}(j-1+\alpha)(j-\alpha)c_{j-1}+2jc_j\}
   t_0^{\alpha-j+1}t^{-\alpha-j-1}\chi_{j-1/2-n/2}(k(\cdot)).
$$
If we choose the coefficients $c_j$ so that
\begin{equation}\label{hd6}
c_0=1\quad\mbox{and}\quad {1\ov2}(j-1+\alpha)(j-\alpha)c_{j-1}+2jc_j=0,\ j\geq 1,
\end{equation}
the then (\ref{hd5}) reduces to
\begin{equation}\label{hd5'}
{\cal P}_\alpha E_\alpha=\delta(x,t-t_0),
\end{equation}
that is $E_\alpha$ is a fundamental solution of $\cal P_\alpha.$
Now recalling notations
\begin{equation}\label{hd6'}
(a)_0=1,\quad (a)_j=a(a+1)\cdots(a+j-1)={\Gam(a+j)\ov\Gam(a)},\ j\geq 1,
\end{equation} 
it follows from (\ref{hd6}) that 
\begin{equation}\label{hd6''}
c_j=(-{1\ov4})^j{(\alpha)_j(1-\alpha)_j\ov j!},\ j\geq 0.
\end{equation}
Hence we may rewrite (\ref{hd4}) as 
\begin{equation}\label{hd4'}
E_\alpha(x,t;0,t_0)=\pi^{1/2-n/2}({t_0\ov t})^\alpha\Phi_\alpha(x,t),
\end{equation}
where 
\begin{equation}\label{hd7}
\Phi_\alpha(x,t)=\sum_{j=0}^\infty {(\alpha)_j(1-\alpha)_j\ov j!}(-{1\ov4t_0t})^j\chi_{j+1/2-n/2}(k(\cdot)).
\end{equation}
This series converges for $|k(\cdot)/4t_0t|<1.$ $\Phi_\alpha(x,t)$ is the {\em hypergeometric
distribution} introduced by Delache and Leray in \cite{dl}. Hypergeometric distributions were also considered by
Gelfand and Shilov in \cite{gs}. 

The expression of $\Phi_\alpha$ depends on the space dimension $n.$ To see this consider three cases.  

\msk
{\bf Case I:} $n=1.$ We have 
\begin{equation}\label{hd8}
\Phi_\alpha(x,t)=\chi_0(k(\cdot)) +
 \sum_{j=1}^\infty {(\alpha)_j(1-\alpha)_j\ov j!}(-{1\ov4t_0t})^j\chi_{j}(k(\cdot)).
\end{equation}
From Euler's formula (\ref{gn8}) it follows that $\chi_j(s)=s^j\chi_0(s)/j!,\ j\geq 0.$ By recalling the expression
of $\chi_{j}(k(\cdot))$ we rewrite (\ref{hd8}) as follows
\begin{eqnarray}\label{hd8'}
\Phi_\alpha(x,t) &=&\chi_0(k(\cdot))
 \sum_{j=0}^\infty {(\alpha)_j(1-\alpha)_j\ov j!j!}({(t-t_0)^2-x^2\ov-4t_0t})^j\\
  &=&  \chi_0(k(\cdot))F(\alpha,1-\alpha,1;{(t-t_0)^2-x^2\ov-4t_0t}) . \nonumber
\end{eqnarray}
 
\msk
{\bf Case II:} $n$ even. We have
\begin{equation}\label{hd9}
\Phi_\alpha(x,t)=\chi_{1/2-n/2}(k(\cdot)) +
  \sum_{j=1}^\infty {(\alpha)_j(1-\alpha)_j\ov j!}(-{1\ov4t_0t})^j\chi_{j+1/2-n/2}(k(\cdot)).
\end{equation}
From Euler's formula (\ref{gn8}) it follows by induction that 
$$
s^j\chi_q(s)=(q+1)_j\chi_{q+j}(s), 
$$
for all integer $j\geq 0.$ Inserting the corresponding formula with $q=1/2-n/2$ into (\ref{hd9}) we obtain
\begin{eqnarray}\label{hd9'}
\Phi_\alpha(x,t) &=&\chi_{1/2-n/2}(k(\cdot))
 \sum_{j=0}^\infty {(\alpha)_j(1-\alpha)_j\ov (3/2-n/2)_j\,j!}({(t-t_0)^2-|x|^2\ov-4t_0t})^j\\
  &=&  \chi_{1/2-n/2}(k(\cdot))F(\alpha,1-\alpha,{3\ov2}-{n\ov2};{(t-t_0)^2-|x|^2\ov-4t_0t}) . \nonumber
\end{eqnarray}

\msk
{\bf Case III:} $n$ odd $>1.$ Let $n=2m+1, m\geq 1.$ Note that in this case $1/2-n/2=-m,$ a negative integer. We
split $\Phi_\alpha$ into two terms:
\begin{eqnarray}\label{hd10}
\Phi_\alpha(x,t)
   &=&\sum_{j=0}^{m-1} {(\alpha)_j(1-\alpha)_j\ov j!}(-{1\ov4t_0t})^j\chi_{j-m}(k(\cdot)) +\\
   &+&\sum_{j=m}^\infty {(\alpha)_j(1-\alpha)_j\ov j!}(-{1\ov4t_0t})^j\chi_{j-m}(k(\cdot)). \nonumber
\end{eqnarray}
Whenever $j-m<0,$ $\chi_{j-m}(k(\cdot))=\delta^{(m-j-1)}(k(\cdot))$ is a distribution concentrated
on the surface of the semi-cone $C.$ Thus the first term in (\ref{hd10}) corresponds to a finite sum of
distributions supported by the boundary of $C.$ 

Recalling that $\chi_{j}(s)=s^{j}\chi_0(s)/j!$ and setting $j'=j-m,$ rewrite the second term in (\ref{hd10}) as
$$
S=\chi_0(k(\cdot))({1\ov-4t_0t})^m\sum_{j'=0}^\infty {(\alpha)_{j'+m}(1-\alpha)_{j'+m}\ov( j'+m)!\,j'!}
    ({k(\cdot)\ov-4t_0t})^{j'}
$$
Now 
$(\alpha)_{j'+m}=(\alpha)_m(\alpha+m)_{j'},$ $(1-\alpha)_{j'+m}=(1-\alpha)_m(1-\alpha+m)_{j'},$
and $(j'+m)!=m!(m+1)_{j'}.$ Therefore
$$
S=\chi_0(k(\cdot))c_m({1\ov-4t_0t})^mF(\alpha+m,1-\alpha+m,m+1,
    {(t-t_0)^2-|x|^2\ov-4t_0t}),
$$
where $c_m=(\alpha)_m(1-\alpha)_m/m!.$ Thus the expression (\ref{hd10}) for $\Phi_\alpha$ becomes
\begin{eqnarray}\label{hd10'}
\Phi_\alpha(x,t)
   &=&\sum_{j=0}^{m-1} {(\alpha)_j(1-\alpha)_j\ov j!}(-{1\ov4t_0t})^j\delta^{(m-j-1)}(k(\cdot)) +\\
   &+&\chi_0(k(\cdot))c_m({1\ov-4t_0t})^mF(\alpha+m,1-\alpha+m,m+1,
    {(t-t_0)^2-|x|^2\ov-4t_0t}). \nonumber
\end{eqnarray}

{\bf Remarks} 1) The support of all fundamental solutions above described is the closure of the
semi-cone $C$ defined at the end of Section \ref{gn}. In the case $n$ odd integer $>1,$ besides the term
that contains the hypergeometric function whose support is the closure of $C$ there are a finite number of terms
whose support is the boundary of $C.$  

2) Formula (\ref{hd10'}) can be viewed as a derivative with respect to $k(\cdot)$ of a certain hypergeometric
distribution. More precisely, consider the hypergeometric distribution $\chi_0(s)F(a,b,1;rs)$ where r is a real or
complex parameter. The following formula holds
\begin{eqnarray}\label{hd11}
{d^m\ov ds^m}[\chi_0(s)F(a,b,1;rs)]&=& \sum_{j=0}^{m-1}{(a)_j(b)_j\ov j!}r^j\delta^{(m-j-1)}(s) +\\
     &+& \chi_0(s)c_mr^mF(a+m,b+m,m+1;rs).   \nonumber
\end{eqnarray}

Indeed, just note that if $f(s)$ is a smooth function defined near $s=0,$ then $f(s)\delta(s)=f(0)\delta(s),$ and
whenever $c\neq 0,-1,-2,\cdots$ one has
$$
{d\ov dz}F(a,b,c;z)={ab\ov c}F(a+1,b+1,c+1;z).
$$

Thus we may rewrite (\ref{hd10'}) as a derivative:
\begin{equation}\label{hd12}
\Phi_\alpha(x,t)={d^m\ov d(k(\cdot))^m}[\chi_0(k(\cdot))F(\alpha,1-\alpha,1,{k(\cdot)\ov -4t_0t})].
\end{equation}
Formulas (\ref{hd11}) and (\ref{hd12}) are analogous to formulas considered by Gelfand and
Shilov in \cite{gs} and involving complex order derivatives of hypergeometric distributions of the type      
$\chi_0(s)F(a,b,c;s).$
\msk

Returning to the operator ${\cal T}_h$ formula (\ref{hd4'})  with $\alpha=1/6$ gives us a fundamental solution 
relative to the point $(0,t_0):$
\begin{equation}\label{hd13}
E_{1/6}(x,t;0,t_0)=\pi^{1/2-n/2}({t_0\ov t})^{1/6}\Phi_{1/6}(x,t).
\end{equation}
In view of our remarks at the end of Section 1, the distribution 
$$
(t/t_0)^{1/3}E_{1/6}(x,t;0,t_0)=\pi^{1/2-n/2}({t\ov t_0})^{1/6}\Phi_{1/6}(x,t)
$$
is then a fundamental solution of ${\cal T}_h^*$ relative to the same point. Motivated by formula (\ref{hd1}) 
we define the distribution $E^\shp$ by
$$
2({3\ov2})^{1/3}t^{1/3}E^\shp(x,t;0,t_0)=\pi^{1/2-n/2}({t\ov t_0})^{1/6}\Phi_{1/6}(x,t),
$$
or,
\begin{equation}\label{hd14}
E^\shp(x,t;0,t_0)={\pi^{1/2-n/2}\ov2^{1/3}3^{1/3}}({1\ov4t_0t})^{1/6}\Phi_{1/6}(x,t).
\end{equation}

In the next sections, we derive from this formula fundamental solutions to the Tricomi operator (\ref{gn1}) and
relative to a point $(0,b), b<0.$ We must distinguish three cases: (I) $n=1$ which corresponds to the classical
Tricomi operator, (II) $n$ an {\em even}, and (III) $n$ {\em odd} $>1.$ In order to simplify notations we write, in
what follows, $E(x,t;0,t_0)$ instead of $E^\shp(x,t;0,t_0).$

\section{The classical Tricomi operator}\setcounter{equation}{0}\label{cl}

If $n=1,$ then (\ref{gn1}) is the classical Tricomi operator in two variables. For this operator we will 
obtain two distinct fundamental solutions: one with support in a region entirely contained in the hyperbolic
half-plane and the other with support in the complement of that region. From formula (\ref{hd14}) we get
\begin{equation}\label{cl1}
E(x,t;0,t_0)={1\ov2^{1/3}3^{1/3}}({1\ov4t_0t})^{1/6}\Phi_{1/6}(x,t).
\end{equation}
On the other hand, since $F(a,b,c;z)=F(b,a,c;z),$ we get from (\ref{hd8'}) that 
\begin{equation}\label{cl1'}
\Phi_{1/6}(x,t)=\chi_0(k(\cdot))F({5\ov6},{1\ov6},1;-{k(\cdot)\ov4t_0t}).
\end{equation}
Note that  $\chi_0(k(\cdot))$ is the characteristic function of the semi-cone $C.$ Recall that
\begin{equation}\label{cl2}
F(a,b,c;z)=(1-z)^{-b}F(c-a,b,c;{z\ov z-1}). 
\end{equation}
If we set $z=(t-t_0)^2-x^2/-4t_0t,$ then
\begin{equation}\label{cl2a}
1-z={(t+t_0)^2-x^2\ov4t_0t}\quad\mbox{and}\quad
  {z\ov z-1}={(t-t_0)^2-x^2\ov(t+t_0)^2-x^2},
\end{equation}
hence
\begin{equation}\label{cl2b}
F({5\ov6},{1\ov6},1;-{k(\cdot)\ov4t_0t})=({(t+t_0)^2-x^2\ov4t_0t})^{-1/6}
F({1\ov6},{1\ov6},1;{(t-t_0)^2-x^2\ov(t+t_0)^2-x^2}),
\end{equation}
and we rewrite (\ref{cl1}) as follows
\begin{equation}\label{cl3}
E(x,t;0,t_0)=\chi_0(k(\cdot)){((t+t_0)^2-x^2)^{-1/6}\ov2^{1/3}3^{1/3}}
  F({1\ov6},{1\ov6},1;{(t-t_0)^2-x^2\ov(t+t_0)^2-x^2}),
\end{equation}
which is, as we pointed out at the end of Section \ref{hd}, a fundamental solution to ${\cal T}_h^*.$ Since 
$\chi_0(k(\cdot))$ is the characteristic function of the semi-cone $C,$ it follows that $E(x,t;0,t_0)$ is
supported by the closure of $C.$ Moreover, the last factor in formula (\ref{cl3}) represents the hypergeometric
series, because the absolute value of its argument (denoted by $z/z-1$ in formula (\ref{cl2a})) is less than $1.$

If one translates formula (\ref{cl3}) in terms of the variables $x$ and $y$, one obtains a fundamental solution of
the classical Tricomi operator, relative to the point $(0,b), b<0,$ and supported by the closure of the region in
$\mathbb R^2$ that corresponds to the semi-cone $C.$ More specifically, consider the change of variables
(\ref{gn1}) and let $a>0$ be such that $t_0=2(-b)^{3/2}/3=a.$ Then, we have
\begin{equation}\label{cl4}
(t-t_0)^2-x^2=-{1\ov9}[9(x^2-a^2)+12a(-y)^{3/2}+4y^3]
\end{equation}
and
\begin{equation}\label{cl4a}
(t+t_0)^2-x^2=-{1\ov9}[9(x^2-a^2)-12a(-y)^{3/2}+4y^3].
\end{equation}
In what follows and in order to simplify notations, we set
\begin{equation}\label{cl4b}
u=9(x^2-a^2)+12a(-y)^{3/2}+4y^3,\quad v=9(x^2-a^2)-12a(-y)^{3/2}+4y^3.
\end{equation}
One can see that 
$$
u=[3(x-a)+2(-y)^{3/2}][3(x+a)-2(-y)^{3/2}]
$$
where the curve $3(x-a)+2(-y)^{3/2}=0$ is one of the characteristics of $\cal T$ through $(0,b)$ and 
$3(x+a)-2(-y)^{3/2}=0,$ the other. Similarly,
$$
v=[3(x-a)-2(-y)^{3/2}][3(x+a)+2(-y)^{3/2}].
$$
The curve $r_a$ of equation $3(x-a)-2(-y)^{3/2}=0$ corresponds to one of the branches of the characteristic
curve originating from $(a,0)$ while $r_{-a},$ the curve of equation $3(x+a)+2(-y)^{3/2}=0,$ corresponds to one of
the branches of the characteristic originating from $(-a,0).$

It is a matter of verification that the semi-cone $C$ corresponds in $\mathbb R^2$ to the region 
\begin{equation}\label{cl6}
D_{b,-}=\{(x,y)\in\mathbb R_-^2: 9(x^2-a^2)+12a(-y)^{3/2}+4y^3<0, y<b\},
\end{equation}
denoted by $D_I$ in the article \cite{bg1}. One may now represent $E(x,t;0,t_0)$ in terms of $x$ and $y,$ via
the expressions $u$ and $v,$ by
\begin{equation}\label{cl6a}
E_-(x,y;0,b)={\chi_{}}_{D_{b,-}}(x,y)\cdot\frac{(-v)^{-1/6}}{2^{1/3}}
F({1\ov6},{1\ov6},1;{u\ov v})
\end{equation}
where ${\chi_{}}_{D_{b,-}}$ is the characteristic function of $D_{b,-}.$

In order to get another fundamental solution supported by the closure of the complement of $D_{b,-}$ we
introduce, as explained in the Appendix, $\td F(1/6,1/6,1;\zeta),$ the {\em principal branch} of the analytic
continuation of the corresponding hypergeometric series, and define in the whole of $\mathbb R^2$ the function
\begin{equation}\label{cl5}
\td E(x,y;0,b)=\frac{(-v)^{-1/6}}{2^{1/3}}
\td F({1\ov6},{1\ov6},1;{u\ov v}).
\end{equation}
We will see in the Appendix that $\td E(x,y;0,b)$ is locally integrable in $\mathbb R^2,$ singular when $v=0,$
real analytic in $\mathbb R^2\backslash(r_a\cup r_{-a}),$ and a solution ${\cal T}u=0$ in the sense of
distributions. We have the following result:
\begin{theorem}\label{th1}
The distribution $E_-$ defined by
\begin{equation}\label{cl7}
E_-(x,y;0,b) = \left\{\begin{array}{ll}
         \td E(x,y;0,b) & \mbox{in $D_{b,-}$}\\
  \\
                          0 & \mbox{elsewhere}
\end{array} \right .
\end{equation}
is a fundamental solution of the Tricomi  operator ${\cal T}$ relative to the point $(0,b).$ Its support is the
closure of $D_{b,-}$.
\end{theorem}

{\bf Proof.} $E_-$ is just another way of writting the expression (\ref{cl6a}). \hfill$\Box$ 
\msk

{\bf Remarks.} 1. Since in $D_{b,-}$ both $u$ and $v$ are $<0,$  it follows that $E_-(x,y;0,b)$ is real valued. 

2. In \cite{bg1} this fundamental solution was obtained by a method different than the one here described and
based upon the existence of the Riemann function for the reduced hyperbolic operator ${\cal T}_h.$ 

3. $E_-$ is the unique fundamental solution of $\cal T,$ relative to $(0,b)$ whose support is $\bar D_{b,-}.$ 
Indeed, any other such fundamental solution is of the form $E_- +f,$ with ${\cal T}f=0$ and $y\leq b$ on
$\mbox{supp}\, f.$ Since the convolution $E_-*f$ is well defined because the map
$$
\mbox{supp}\,E_-\times\mbox{supp}\,f\owns ((x,y),(x',y'))\to(x+x',y+y') 
$$
is proper, we have
$$
f={\cal T}E_-*f=E_-*{\cal T}f=0.
$$
\msk

As a consequence of Theorem \ref{th1} we obtain one of the fundamental solutions described in \cite{bg}. 
\begin{cor}\label{co1}
As $(0,b)\to(0,0),$ the fundamental solution (\ref{cl7}) converges, in the sense of distributions, to the
fundamental solution
\begin{equation}\label{cl8}
F_-(x,y) = \left\{\begin{array}{ll}
         {\ds{1\over2^{1/3}}}F({1\ov6}, {1\ov6},1;1)|9x^2+4y^3|^{-1/6} & \mbox{in $D_-$}\\
  \\
                          0 & \mbox{elsewhere,}
\end{array} \right .
\end{equation}
where $D_-=\{(x,y)\in\mathbb R^2: 9x^2+4y^3<0\}.$
\end{cor}

Since ${\cal T}\td E=0$ in the sense of distributions, it follows that the distribution $E_- -\td E,$ identically 
zero in the region $D_{b,-},$ is also a fundamental solution of $\cal T.$ Denote by 
$D_{b,+}$ the complement in $\mathbb R^2$ of $D_{b,-}$ and  define the distribution 
\begin{equation}\label{cl9}
E_+(x,y;0,b) = \left\{\begin{array}{ll}
         -\td E(x,y;0,b) & \mbox{in $D_{b,+}$}\\
  \\
                          0 & \mbox{elsewhere.}
\end{array} \right .
\end{equation}
We clearly have
\begin{theorem}\label{th2}
$E_+$ is a fundamental solution of $\cal T$ relative to $(0,b)$ whose support is the closure of the region 
$D_{b,+}.$
\end{theorem}

This fundamental solution is not unique. If we replace the exponential factor in (\ref{cl5}) by $e^{-i\pi/6}$ we
obtain another fundamental solution. Moreover, it does not follow as in the case of $E_-,$ that $E_+$ converges,
as $b\to 0,$ to the fundamental solution $F_+(x,y)$ described in \cite{bg}. In order to obtain such a result, one
needs to consider a suitable linear combination of these two fundamental solutions before taking limits (see
\cite{bg1}).  

We thus have 
\begin{cor}\label{co2}
As $(0,b)\to(0,0),$ a suitable linear combination of fundamental solutions of the type $E_+$ converges, in the
sense of distributions, to the fundamental solution
\begin{equation}\label{cl10}
F_+(x,y) = \left\{\begin{array}{ll}
         {\ds-{1\over2^{1/3}3^{1/2}}}F({1\ov6}, {1\ov6},1;1)(9x^2+4y^3)^{-1/6} & \mbox{in $D_+$}\\
  \\
                          0 & \mbox{elsewhere,}
\end{array} \right .
\end{equation}
where $D_+=\{(x,y)\in\mathbb R^2: 9x^2+4y^3>0\}.$
\end{cor}

\section{The Tricomi operator, $n$ even}\setcounter{equation}{0}\label{ev}

We begin with formula (\ref{hd14}) 
$$
E(x,t;0,t_0)={\pi^{1/2-n/2}\ov2^{1/3}3^{1/3}}({1\ov4t_0t})^{1/6}\Phi_{1/6}(x,t),
$$
where $\Phi_{1/6},$ given by (\ref{hd9'}), is 
\begin{equation}\label{ev1}
\Phi_{1/6}(x,t)=\chi_{1/2-n/2}(k(\cdot))F({5\ov6},{1\ov6},{3\ov2}-{n\ov2};-{k(\cdot)\ov4t_0t}).
\end{equation}
Recalling formulas (\ref{gn7}), (\ref{cl2}) and (\ref{cl2a}) we obtain
 \begin{eqnarray}\label{ev2}
\hskip 1cm\lefteqn{E(x,t;0,t_0)=}\\
&=&{\pi^{1/2-n/2}\ov2^{1/3}3^{1/3}}\,\chi_{1/2-n/2}(k(\cdot))\,
     [(t+t_0)^2-|x|^2]^{-1/6}\times\nonumber \\ 
&\times&  F({2\ov3}-{n\ov2},{1\ov6},{3\ov2}-{n\ov2};{(t-t_0)^2-|x|^2\ov(t+t_0)^2-|x|^2}).\nonumber
\end{eqnarray}

To obtain a fundamental solution to $\cal T$ we represent (\ref{ev2}) in terms of $x$ and $y.$ If we set
\begin{equation}\label{ev3}
u=9(|x|^2\!-\!a^2)\!+\!12a(-y)^{3/2}\!+\!4y^3,\ \ v=9(|x|^2\!-\!a^2)\!-\!12a(-y)^{3/2}\!+\!4y^3,
\end{equation}
then
\begin{equation}\label{ev3a}
(t-t_0)^2-|x|^2=-{1\ov9}u\quad\mbox{and}\quad (t+t_0)^2-|x|^2=-{1\ov9}v.
\end{equation}
These two formulas are the counterpart to (\ref{cl4}) and (\ref{cl4a}) in the case $n=1.$ 

Define, as we did in Section \ref{cl} case $n=1,$ the region
\begin{equation}\label{ev5}
D^n_{b,-}=\{(x,y)\in\mathbb R^{n+1}:9(|x|^2-a^2)+12a(-y)^{3/2}+4y^3<0, y<b\}
\end{equation}
which corresponds to the semi-cone $C,$ and let $\chi_{D^n_{b,-}}$ be its characteristic function.

In terms of $x$ and $y$ the distribution (\ref{ev2})becomes
\begin{equation}\label{ev4}
E_-(x,y;0,b)=c(n)\chi_{D^n_{b,-}}(x,y)(-u)^{1/2-n/2}(-v)^{-1/6}
   F({2\ov3}-{n\ov2},{1\ov6},{3\ov2}-{n\ov2};{u\ov v}),
\end{equation}
where
\begin{equation}\label{ev4'}
c(n)={\pi^{1/2-n/2}\ov2^{1/3}3^{1-n}\Gam({3\ov2}-{n\ov2})}.
\end{equation}

Thus we obtain the following result:
\begin{theorem}\label{tev1} $E_-(x,y;0,b)$ is a fundamental solution of $\cal T$ relative to $(0,b)$ whose support 
is the closure of the region $D^n_{b,-}.$
\end{theorem}

If we let $b\to 0,$ we obtain a fundamental solution of $\cal T$ relative to the origin, namely
\begin{cor}\label{cev1} 
The limit, in the sense of distributions, of $E_-(x,y;0,b)$ as $(0,b)\to(0,0)$ is 
 \begin{eqnarray}\label{ev7}
\lefteqn{F_-(x,y)=}\\
&=&\left\{\begin{array}{ll}
\ds{\pi^{1/2-n/2}\ov2^{1/3}3^{1-n}\Gamma({3\ov2}-{n\ov2})}F({2\ov3}-{n\ov2},{1\ov6},{3\ov2}-{n\ov2};1)
      \,|9|x|^2+4y^3|^{{1\ov3}-{n\ov2}}
       & \mbox{in $D_-^n$}\\
\\
 \hskip 1cm 0 & \mbox{elsewhere,} 
\end{array} \right.\nonumber
\end{eqnarray}
a fundamental solution of ${\cal T}$ relative to the origin whose support is the closure of the region
$D_-^n=\{(x,y)\in{\mathbb R}^{n+1}: 9|x|^2+4y^3<0\}.$ 
\end{cor}

The fundamental solution given by formula (\ref{ev7}) coincides with the fundamental solution given by 
formula (4.2) in Theorem 4.1 of \cite{bc}. The only apparent discrepancy between these 
two formulas is the multiplying constants. In (\ref{ev7}), the multiplicative constant is
\begin{equation}\label{ev8}
A={\pi^{1/2-n/2}\ov2^{1/3}3^{1-n}\Gam({3\ov2}-{n\ov2})}F({2\ov3}-{n\ov2},{1\ov6},{3\ov2}-{n\ov2};1)
\end{equation}
while in \cite{bc}, page 490, the multiplicative constant for $F_-(x,y)$ is 
\begin{equation}\label{ev8a}
C_-={3^n\Gam(4/3)\ov2^{2/3}\pi^{n/2}\Gam({4\ov3}-{n\ov2})}.
\end{equation}
Since 
$$
F(a,b,c;1)={\Gam(c)\Gam(c-a-b)\ov\Gam(c-a)\Gam(c-b)}
$$
a formula that holds whenever  $\mbox{Re}\,c>\mbox{Re}\,b + \mbox{Re}\,a,$ we may rewrite $A$ as
$$
A={\pi^{1/2-n/2}\Gam(2/3)\ov2^{1/3}3^{1-n}\Gam(5/6)\Gam({4\ov3}-{n\ov2})}.
$$
In order for $A=C_-$ one must have the identity
$$
{2^{1/3}\pi^{1/2}\Gam(2/3)\ov3\Gam(5/6)\Gam(4/3)}=1.
$$
But this is a consequence of the following relations for the Gamma function:
$\Gam(z+1)=z\Gam(z),$ $\Gam(2z)=2^{2z-1}\pi^{-1/2}\Gam(z)\Gam(z+1/2),$ and
$\Gam(z)\Gam(1-z)=\pi\csc(\pi z).$ 
\msk

In \cite{bc} we showed that the distribution 
\begin{equation}\label{ev10}
F_+(x,y) = \left\{\begin{array}{ll}
         C_+(9|x|^2+4y^3)^{1/3-n/2} & \mbox{in $D^n_+$}\\
  \\
                          0 & \mbox{elsewhere,}
\end{array} \right .
\end{equation}
where 
\begin{equation}\label{ev10a}
C_+=-{3^{n-2}\Gam({n\ov2}-{1\ov3})\ov2^{2/3}\pi^{n/2}\Gam(2/3)}
\end{equation}
and
\begin{equation}\label{ev10b}
D^n_+=\{(x,y)\in \mathbb R^{n+1}:9|x|^2+4y^3>0\},
\end{equation}
is a fundamental solution of $\cal T$ supported by the closure of $D^n_+.$
It is a matter of verification that the ratio between the constants $C_+$ and $C_-$ is 
\begin{equation}\label{ev11}
{C_+\ov C_-}=-{1\ov2\sqrt3\sin\pi({n\ov2}-{1\ov3})}.
\end{equation}
Therefore the constant $C_+$ in \cite{bc} can also be represented in terms of the above hypergeometric function.

\section{The Tricomi operator, $n$ odd $>1$}\setcounter{equation}{0}\label{od}
Let $n=2m+1$ with $m\geq1.$ Again from formula (\ref{hd14}) we have
\begin{equation}\label{od1}
E(x,t;0,t_0)=A_m({1\ov4t_0t})^{1/6}\Phi_{1/6}(x,t),
\end{equation}
where $A_m=1/2^{1/3}3^{1/3}\pi^m.$ From formula (\ref{hd10'}) $\Phi_{1/6}$ is given by 
\begin{eqnarray}\label{od2}
\Phi_{1/6}(x,t)
   &=&\sum_{j=0}^{m-1} c_j(-{1\ov4t_0t})^j\delta^{(m-j-1)}(k(\cdot)) +\\
   &+&\chi_0(k(\cdot))c_m({1\ov-4t_0t})^mF(m+{5\ov6},m+{1\ov6},m+1,
    {k(\cdot)\ov-4t_0t}). \nonumber
\end{eqnarray}
with
\begin{equation}\label{od2'}
c_j={\Gam(j+5/6)\Gam(j+1/6)\ov\Gam(5/6)\Gam(1/6)\Gam(j+1)},\quad 0\leq j\leq m.
\end{equation}
In view of (\ref{cl2}) and (\ref{cl2a}), the hypergeometric function in (\ref{od2}) is equal to
$$
({t+t_0)^2-|x|^2\ov4t_0t})^{-m-1/6}
   F({1\ov6},m+{1\ov6},m+1;{(t-t_0)^2-|x|^2\ov(t+t_0)^2-|x|^2})
$$
and we may rewrite $E(x,t;0,t_0)$ as 
\begin{equation}\label{od3}
 E(x,t;0,t_0)=
\end{equation}
$$
  =A_m\sum_{j=0}^{m-1}(-1)^j c_j(4t_0t)^{-j-1/6}\delta^{(m-j-1)}(k(\cdot)) +
$$
$$
  +(-1)^mA_mc_m((t+t_0)^2-|x|^2)^{-m-1/6}F({1\ov6},m+{1\ov6},m+1,
    {(t-t_0)^2-|x|^2\ov(t+t_0)^2-|x|^2})\chi_0(k(\cdot)). \nonumber
$$
Note that all terms in the sum contain distributions of the form $\delta^{(q)}(k(\cdot))$ which 
are supported by the surface of the semi-cone $C.$ However, the support of the last term in (\ref{od3}) is the
closure of $C.$

In \cite{gs} Gelfand and Shilov introduced the distribution $\delta(P)$ supported by the surface $S$ given by$P=0,$
where $P$ is a smooth function such that $\nabla P\neq 0$ on $S.$ In particular, they proved that if $a(\cdot)$ is a
nonvanishing function, then
\begin{equation}\label{od4}
\delta^{(q)}(aP)=a^{-(q+1)}\delta^{(q)}(P).
\end{equation}
These results extend to our case, where $P=k(\cdot)$ has a singular point at $(0,t_0).$ We have the following
\begin{lemma}\label{odl}
For all $0\leq j\leq m-1,$
$$
(4t_0t)^{-j-1/6}\delta^{(m-j-1)}(k(\cdot))=
    (4t_0t)^{5/6}\delta^{(m-j-1)}((4t_0t)^{j+1/m-j}k(\cdot)).
$$
\end{lemma}

{\bf Proof.} Indeed we have
$$
(4t_0t)^{-j-1/6}\delta^{(m-j-1)}(k(\cdot))=(4t_0t)^{5/6}(4t_0t)^{-(j+1)}\delta^{(m-j-1)}(k(\cdot))=
$$
$$
=\!(4t_0t)^{5/6}[(4t_0t)^{j+1/m-j}]^{-(m-j)}\delta^{(m-j-1)}(k(\cdot))
  \!=\!(4t_0t)^{5/6}\delta^{(m-j-1)}((4t_0t)^{j+1/m-j}k(\cdot)),
$$
by virtue of (\ref{od4}) and the fact that $4t_0t\neq0$ in the region $t-t_0>|x|.$\hfill$\Box$
\msk

As a consequence of this lemma, all terms that contain derivatives of $\delta$ in (\ref{od3}) tend to zero, as
$t_0\to0.$ By taking limits, it follows that the distribution
\begin{equation}\label{od5}
E(x,t;0,0)=
\end{equation}
$$
={(-1)^m\ov2^{1/3}3^{1/3}\pi^m}{\Gam(m+5/6)\Gam(m+1/6)\ov\Gam(5/6)\Gam(1/6)\Gam(m+1)}
 F({1\ov6},m+{1\ov6},m+1;1)(t^2-|x|^2)^{-m-1/6}
$$
is a fundamental solution of ${\cal T}_h$ supported by the closure of the semi-cone $\{(x,t):t>|x|\}.$

As we did in the previous sections, we rewrite $E(x,t;0,t_0)$ in terms of the variables $x$ and $y.$ From
formulas (\ref{ev3}) and (\ref{ev3a}) we derive that $4t_0t=(v-u)/9$ and, following Gelfand and Shilov's
notations, we replace $\delta^{(q)}(k(\cdot))$ by $\delta^{(q)}(u(\cdot)),$ with the understanding that
$u(\cdot)$ now means $u(x,y),$ with $y\leq b.$ Thus (\ref{od3}) becomes
\begin{equation}\label{od3a}
 E_-(x,y;0,b)=
\end{equation}
$$
  =A_m\sum_{j=0}^{m-1}(-1)^j c_j(\frac{v-u}{9})^{-j-1/6}\delta^{(m-j-1)}(u(\cdot)) +
$$
$$
  +(-1)^mA_mc_m(-\frac{v}{9})^{-m-1/6}F({1\ov6},m+{1\ov6},m+1,\frac{u}{v})\chi_{D^n_{b,-}}(x,y) \nonumber
$$
where $\chi_{D^n_{b,-}}$ is the characteristic function of the set (\ref{ev5}).Then the following result holds:

\begin{theorem}\label{odt}
The distribution $E_-(x,y;0,b)$ is a fundamental solution of $\cal T$ relative to $(0,b)$ 
supported by the closure of the set $D^n_{b,-}.$ 
\end{theorem}
\msk

Note that in (\ref{od3a}) all terms inside the summation are supported by the boundary of $D^n_{b,-}$ 
while the last   term is supported by the closure of $D^n_{b,-}.$ If we let $b\to0,$ we obtain at the limit, 
the fundamental solution  $F_-(x,y)$ described in our previous paper \cite{bc}, namely

\begin{theorem}\label{odt1} The distribution
\begin{equation}\label{od5a}
F_-(x,y) = \left\{\begin{array}{ll}
\ds{3^n\Gamma({4\ov3})\ov2^{2/3}\pi^{n/2}\Gamma({4\ov3}-{n\ov2})}
      \,|9|x|^2+4y^3|^{{1\ov3}-{n\ov2}}
       & \mbox{in $D_-^n$}\\
\\
 \hskip 1cm 0 & \mbox{elsewhere,} 
\end{array} \right.
\end{equation}
supported by the closure of the region $D_-^n=\{(x,y)\in{\mathbb R}^{n+1}: 9|x|^2+4y^3<0\},$ is a
fundamental solution of $\cal T.$
\end{theorem}

{\bf Proof.} Recall that $t=2(-y)^{3/2}/3$ and that $t^2-|x|^2={1\ov9}(-9|x|^2-4y^3).$ Hence, the
right hand-side of (\ref{od5}) equals
$$
A F({1\ov6},m+{1\ov6},m+1;1)|9\,|x|^2+4y^3|^{-m-1/6}
$$
where
\begin{equation}\label{od6}
A={(-1)^m3^{2m}\ov2^{1/3}\pi^m}{\Gam(m+5/6)\Gam(m+1/6)\ov\Gam(5/6)\Gam(1/6)\Gam(m+1)}.
\end{equation} 
Now the exponent $-m-1/6$ equals $1/3-n/2$ because $n=2m+1.$ On the other hand, it is a matter of
verification that the constant 
\begin{equation}\label{od7}
{3^n\Gamma({4\ov3})\ov2^{2/3}\pi^{n/2}\Gamma({4\ov3}-{n\ov2})}
\end{equation}
(denoted by $C_-$ in \cite{bc}) which appears in (\ref{od5a}) is the same as $A.$\hfill$\Box$

\section{Appendix}\setcounter{equation}{0}\label{app}
We are going to prove that the function $\td E(x,y;0,b)$ defined by formula (\ref{cl5}) in Section \ref{cl} is locally
integrable in $\mathbb R^2,$ singular when $v=0,$ real analytic in $\mathbb R^2\backslash(r_a\cup r_{-a}),$ and a
solution of ${\cal T}w=0$ in the sense of distributions. Recall that $r_a$ is the characteristic curve
$3(x-a)-2(-y)^{3/2}=0$ originating from $(a,0),$ and $r_{-a},$ the characteristic curve $3(x+a)+2(-y)^{3/2}=0,$
originating from
$(-a,0).$

Following Whittaker and Watson \cite{ww}, let $\alpha,$ $\beta,$ and $\gam$ be complex numbers,
$\gam\neq0,-1,-2,\cdots,$ and let
\begin{equation}\label{app2}
(\alpha)_0=1, (\alpha)_n=\alpha(\alpha+1)\cdots(\alpha+n-1)=\frac{\Gam(\alpha+n)}{\Gam(\alpha)}.
\end{equation}
The power series
\begin{equation}\label{app2a}
F(\alpha,\beta,\gam;\zeta)=\sum_{n=0}^\infty\frac{(\alpha)_n(\beta)_n}{(\gam)_n n!}\zeta^n
\end{equation}
is called the {\em hypergeometric series}. The ratio test guarantees absolute convergence for $|\zeta|<1.$ If 
$\Re(\gam-\alpha-\beta>0,$ then the series converges for $|\zeta|\leq1$ and 
\begin{equation}\label{app2b}
F(\alpha,\beta,\gam;1)=\frac{\Gam(\gam)\Gam(\gam-\alpha-\beta)}{\Gam(\gam-\alpha)\Gam(\gam-\beta)}.
\end{equation}

{\em Barnes' contour integral} defines a single-valued analytic function of $\zeta$ in the region
$|\arg(-\zeta)|<\pi,$ that is, $\mathbb C$ minus the positive real axis, which gives the {\em principal branch} of 
the analytic continuation of the hypergeometric series $F(\alpha,\beta,\gam;\zeta).$ More precisely we quote the
following theorem whose proof is found in \cite{ww}.
\begin{theorem}\label{bar}
{\bf (Barnes)} The integral
\begin{equation}\label{app3}
\frac{1}{2\pi i}\int_{-i\infty}^{i\infty}\frac{\Gam(\alpha+s)\Gam(\beta+s)\Gam(-s)}{\Gam(\gam+s)}
   (-\zeta)^s\,ds,
\end{equation}
where the contour of integration is curved (if necessary) to ensure that the poles
of $\Gam(\alpha+s)\Gam(\beta+s),$ i.e., $s=-\alpha-n, -\beta-n, n=0,1,2,\cdots,$ lie on the left of the contour and
the poles of $\Gam(-s),$ i.e., lie on the right of the contour, defines a single-valued analytic function in the region 
$|\arg(-\zeta)|<\pi.$ Moreover, in the unit disk $|\zeta|<1,$ it coincides with the hypergeometric series
$$
\frac{\Gam(\alpha)\Gam(\beta)}{\Gam(\gam)}F(\alpha,\beta,\gam;\zeta).
$$
\end{theorem}

Following traditional practice we use the notation $F(\alpha,\beta,\gam;\zeta)$ to denote either the
hypergeometric series or the principal branch of its analytic continuation, and call it the hypergeometric function. 

Barnes' integral may also be used to obtain a representation of the hypergeometric function in the form of a power
series in $\zeta^{-1},$ convergent when $|\zeta|>1.$ By choosing a suitable contour of integration one can prove
(see \cite{ww}) that if $\alpha-\beta$ is not an integer or zero, then 
\begin{eqnarray}\label{app4}
F(\alpha,\beta,\gam;\zeta)&=&A(-\zeta)^{-\alpha}F(\alpha,1-\gam+\alpha,1-\beta+\alpha;\zeta^{-1})\\
&=&A(-\zeta)^{-\beta}F(\beta,1-\gam+\beta,1-\alpha+\beta;\zeta^{-1}), \nonumber
\end{eqnarray}
where $A$ and $B$ are suitable constants and $|\arg(-\zeta)|<\pi.$ This formula also describes the
asymptotic behaviour of the function $F(\alpha,\beta,\gam;\zeta)$ near $|\zeta|=\infty.$ If
$\alpha-\beta$ is an integer or zero, formula (\ref{app4}) must be modified because some of the poles of
$\Gam(\alpha+s)\Gam(\beta+s)$ are double poles. The reader should find the expression for
$F(\alpha,\beta,\gam;\zeta)$ in \cite{erd}, chapter on hypergeometric functions. In the case that interests
us, that is, $\alpha=\beta,$ that expression is 
\begin{equation}\label{app5}
F(\alpha,\alpha,\gam;\zeta)=(-\zeta)^{-\alpha}[\log(-\zeta)U(\zeta)+V(\zeta)],
\end{equation}
where $|\arg(-\zeta)|<\pi,$ and both $U(\zeta)$ and $V(\zeta)$ are power series in $\zeta^{-1}$ convergent for
$|\zeta|>1.$ The reader expressions are found in \cite{erd} or \cite{bg1}. We also mention that
if $\Re(\gam-\alpha-\beta)>0,$ we have convergence for $|\zeta|\geq1.$

From the above results and in particular from (\ref{app5}) it follows that
$$
\td E(x,y;0,b)=\frac{(-v)^{-1/6}}{2^{1/3}}F(\frac{1}{6},\frac{1}{6}, 1;\frac{u}{v}),
$$
with $u$ and $v$ defined by (\ref{cl4b}), is locally integrable in $\mathbb R^2,$ singular when $v=0,$ and
real analytic in $\mathbb R^2\backslash(r_a\cup r_{-a}).$ 

It remains to prove that $\td E(x,y;0,b)$ is a solution to ${\cal T}w=0$ in the sense of distributions. For this we
need several results proved in the paper \cite{bg}. In that paper we showed that the function
$$
{\cal E}(\ell,m;\ell_0,m_0)=(\ell-m)^{-1/6}(\ell_0-m)^{-1/6}
       F(\frac{1}{6},\frac{1}{6}, 1;\frac{(\ell-\ell_0)(m-m_0)}{(\ell-m_0)(m-\ell_0})
$$
is a classical solution of
$$
{\cal T}_h w=\frac{\pl^2w}{\pl\ell\pl m}-\frac{1/6}{\ell-m}(\frac{\pl w}{\pl\ell}-\frac{\pl w}{\pl m})=0,
$$
the reduced hyperbolic Tricomi equation. Here 
$$
\ell=x+\frac{2}{3}(-y)^{3/2},\qquad m=x-\frac{2}{3}(-y)^{3/2}
$$
are the characteristic coordinates. Now, except for the constant $1/2^{1/3},$ $\td E(x,y;0,b)$ is obtained from
${\cal E}(\ell,m;\ell_0,m_0)$ after replacement of $\ell$ and $m$ by their expressions above and by setting
$\ell_0=-m_0=2(-b)^{3/2}/3.$ Thus, away from the set $\{v=0\}=r_a\cup r_{-a},$ $\td E(x,y;0,b)$ is a classical
solution of ${\cal T}w=0.$

To show that ${\cal T}\td E=0,$ in the sense of distributions, we have to contend with the fact that $\td E(x,y;0,b)$
has logarithmic singularities along the two characteristics $r_{-a}$ and $r_a,$ or, equivalently, that ${\cal
E}(\ell,m;\ell_0,m_0)$ has logarithmic singularities along the lines $\ell=-\ell_0$ and $m=\ell_0.$ Since, as we have
remarked, ${\cal T}\td E=0$ away from the characteristics $r_{-a}$ and $r_a,$, in order to prove that ${\cal T}\td
E=0$ in the sense of distributions, it suffices to prove that
\begin{equation}\label{app6}
\langle\td E,{\cal T}\phi\rangle=\int\int_{\mathbb R^2}\td E{\cal T}\phi dx\,dy=0
\end{equation}
for all $\phi\in\smoc(\mathbb R^2)$ whose support intersects at least one of the characteristics $r_{-a}$ or 
$r_a.$ If supp\,$\phi$ does not intersect either of these characteristics, then (\ref{app6}) is automatically 
satisfied. 

Suppose that supp\,$\phi$ is contained in an open disk $D$ centered, say at $(a,0),$  and with radius $R.$
Let $0<r<R$ and denote by $D_\ep$ the set of points of $D$ at a distance $>\ep$ from the characteristic $r_a.$
Then, from Green's formula for $\cal T$ (see \cite{bg1}, formula (4.5)) one gets
\begin{eqnarray}\label{app7}
\int\int_{D}{\td E}{\cal T}\phi\,dx\,dy &=&
       \lim_{\ep\to0}\int\int_{D_\ep}{\td E}{\cal T}\phi\,dx\,dy\\
  &=& \lim_{\ep\to0}\int_{\Gam_\ep\cup\gam_\ep\cup\Gam'_\ep} {\td E}(y\phi_x\,dy-\phi_y\,dx)-
         \phi(y{\td E}_x\,dy- {\td E}_y\,dx),      
\nonumber
\end{eqnarray}
where $\Gam_\ep$ is the charateristic $3(x-\alpha+\ep)-2(-y)^{3/2}=0,$
$\gam_\ep$ the circumference center at $\alpha$ with radius $\ep,$ and $\Gam'_\ep$ the characteristic
$3(x-\alpha-\ep)-2(-y)^{3/2}=0.$ In order to prove (\ref{app6}) we must prove that the last limit in
(\ref{app7}) is zero. 

Most details of the proof are to be found in Section 4 of the paper \cite{bg1}. We just point out that the integrand 
in (\ref{app7}) remains bounded along $\gam_\ep$ thus, along this contour, the integral tends to zero with $\ep.$
Along both $\Gam_\ep$ and $\Gam'_\ep$ we must take into account the asymptotic behaviour of
$F(1/6,1/6,1;\zeta)$ and its derivative $F(7/6,7/6,2;\zeta),$ at $\zeta=\infty,$ according with (\ref{app5}). 
It turns out that at the limit, the values of these integrals cancel each other and this completes the proof.

\end{document}